\documentclass[12pt]{amsart}
\usepackage{amsmath}
\usepackage{amssymb}
\usepackage{amsfonts}
\usepackage{amsthm}
\usepackage{verbatim}
\usepackage{amscd}
\usepackage{cite}
\usepackage{leftidx}
\usepackage{enumerate}
\usepackage{txfonts}
\usepackage{manfnt}
\usepackage{amscd}
\usepackage[mathscr]{eucal}
\usepackage{hyperref}
\usepackage{datetime2}
\usepackage{mathpazo}
\def\cP{\mathcal P}
\def\cX{\mathcal X}

\newtheorem{thm}{Theorem} 

\newtheorem*{thm*}{Theorem}
\newtheorem*{prop*}{Proposition}
\newtheorem{cor}[thm]{Corollary}
\newtheorem*{cor*}{Corollary}

\newtheorem{lem}[thm]{Lemma}
\newtheorem*{lem*}{Lemma}

\newtheorem*{claim*}{Claim}

\theoremstyle{remark}

\newtheorem{rem}[thm]{Remark}
\newtheorem*{rem*}{Remark}
\newtheorem{crit-rem}[thm]{Critical remark}

\newtheorem{example}[thm]{Example}
\newtheorem*{example*}{Example}

\newtheorem*{defn*}{Definition}
\newtheorem*{con*}{Conjecture}

\def\inv{^{-1}}

\def\cE{\mathcal E}
\def\cC{\mathcal C}

\def\refp #1.{(\ref{#1})}

\newcommand\ceil [1] {\lceil #1 \rceil}

\newcommand{\A}{\mathcal{A}}

\newcommand{\kk}{\mathbf{k}}

\def\sbr #1.{^{[#1]}}
\def\sfl #1.{^{\lfloor #1\rfloor}}

\def\inv{^{-1}}
\def\?{{\bf{??}}}

\def\A{\Bbb A}

\def\C{\mathbb C}
\def\P{\mathbb P}

\def\O{\mathcal O}

\def\rk{\text{rk}}

\def\g{\mathfrak g}

\def\1/2{\frac{1}{2}}

\def\I{\mathcal{ I}}

\def\2{{[2]}}
\def\l{\ell}
\def\nl{\newline}

\def\<{\langle}
\def\>{\rangle}

\def\2{{[2]}}
\def\l{\ell}
\def\Def{\mathrm{Def}}

\def\scl #1.{^{\lceil#1\rceil}}
\def\spr #1.{^{(#1)}}
\def\sbc #1.{^{\{#1\}}}

\def\subpr#1.{_{(#1)}}

\def\beq{\begin{equation*}}
\def\eeq{\end{equation*}}

\newcommand{\mlog}[1]{\langle-\log {#1} \rangle}
\def\g3{{\Gamma\spr 3.}}

\def\Def{\mathrm {Def}}

\newcommand{\eqspl}[2]{
\begin{equation}\label{#1}
\begin{split}
#2\end{split}\end{equation}}

\newcommand{\exseq}[3]{
0\to #1\to #2\to #3\to 0
}

\newcommand{\beginalphaenum}{
\begin{enumerate}\renewcommand{\labelenumi}{ }
\item \begin{enumerate}
}

\def\eex{\end{rm}\end{example}}

\begin{document} 
	\title{Rigid and stably balanced curves\\ on 
	Calabi-Yau and general-type  hypersurfaces}
	\author 
	{Ziv Ran}
%
%
	\thanks{arxiv.org  2209.05410 }
	\date {\DTMnow}
%
%
	\address {\nl UC Math Dept. \nl
	Skye Surge Facility, Aberdeen-Inverness Road
	\nl
	Riverside CA 92521 US\nl 
	ziv.ran @  ucr.edu\nl
	\url{https://profiles.ucr.edu/app/home/profile/zivran}
	}
	
	 \subjclass[2010]{14j70, 14j32}
	\begin{abstract}
	A curve $C$ on a variety $X$ is stably balanced if the slopes of the Harder-Narasimhan filtration
	of its normal bundle $N$
	are contained in an interval of length 1. 
	For a general nonspecial curve of  given genus $g\geq 1$ and large enough degree in $\P^n $,
	and a general hypersurface $X$ of large enough degree containing $C$, we prove that $C$ is
	stably balanced and rigid on $X$
		and that the 
		  family of hypersurfaces $X$ thus obtained  is smooth of the expected codimension
		  $h^1(N)$ in the space of hypersurfaces.
		\end{abstract}
	\maketitle
	\section*{Declarations of Compliance with Ethical Standards}
	\subsection*{Conflict of Interest Statement}
	  The author has no affiliation with any organization with a direct or indirect 
	financial or non-financial interest in the subject matter discussed in the manuscript
	\subsection*{Data Availability Statement}
	There is no data set associated with this paper
	\subsection*{Human Participants and/or Animals}
None were involved in this research.	
	\section*{Introduction}
	\subsection{Known results: interpolation, rigidity, stability}
	A fundamental object of interest in algebraic geometry is a pair $(C, X)$
	consisting of a curve $C$
	 on a hypersurface $X$ 
	in $\P^n$.
	Generally, given a curve $C$ on a variety $X$, 
important attributes of $C\subset X$ such as the motion of $C$ on $X$ 
 reside in in properties of the normal bundle
	$N=N_{C/X}$, i.e. the 1st-order infinitesimal neighborhood of $C$ in $X$
	which controls motions both local and global of $C$ in $X$.
Results on $N$ for curves on general  hypersurfaces have thus far been largely 
restricted to the case where $X$ 
is in the  Fano range, i.e of degree $\leq n$,
so that  $N$ is 'positive' (e.g. $\chi(N)\geq 0$). In this case, going beyond 
the regularity property, i.e. $H^1(N)=0$, 
one natural property of $N$ is interpolation or balance,
	which states that for a general effective divisor $D$ of any degree, 
	one has either $H^0(N(-D))=0$ or $H^1(N(-D))=0$,
	 and there are now
	a number of papers that prove this property when $X$ is $\P^n$ itself 
	(see Atanasov-Larson-Yang \cite{alyang} ) or a Fano
	hypersurface \cite{caudatenormal}, \cite{elliptic}. 
	
Much less is known 
about curves with good normal bundles $N$
when the hypersurface   $X$ is in 
the so-called CYGT range (degree $>n$ in $\P^n$) where $N$ 
need not be positive in any sense (in fact, one has $\chi(N)\leq 0$).
	One important good property of $N$ in case $\chi(N)\leq 0$
	is rigidity, i.e. $H^0(N)=0$. 
	On Calabi-Yau complete intersection 3-folds (CYCI3fs) $X$, where $\chi(N)=0$ for any $g$,
	there are examples in the literature of rigid curves
	going back to Clemens \cite{clemens-quintic}, 
	as extended by Katz \cite{katz-5tic},
	who constructed infinitely many rigid rational curves on the general
	quintic in $\P^4$. This was extended to
	rational curves on all  CICY3fs
	 by Ekedahl, Johnsen and Sommervoll
	\cite{ekedahl1}, and subsequently, based on Clemens's method,
	 by Kley \cite{kley} to  any CICY3f $X$ and curves of
	sufficiently low genus $g$  (depending on $X$, e.g. $g<35$ for the quintic) and sufficiently
	high degree (depending on the genus). For the quintic threefold, any $g\geq 0$ and
	large $d$, Zahariuc \cite{zahariuc} constructs
	a degree- $d$ map from a smooth curve of genus $g$
	to $X$ which is set-theoretically  isolated (it is not proved the map
	is an embedding or  infinitesimally rigid).\par
	For (CY) hypersurfaces of degree $n+1$ in $\P^n, n\geq 4$,
	and some CY complete intersecrtions, we proved in \cite{relreg}
	the existence of rational curves with balanced normal bundle of the form $(n-4)\O\oplus2\O(-1)$.
	\par
	Another property of $N$  which makes sense irrespective of positivity is (semi)stability,
		but results establishing this have so far been restricted to the case $X=\P^n$
		or $X$ a (Fano) hypersurface  in $\P^{n+1}$
		\cite{coskun-larson-vogt-normal}, \cite{semistable}.\par
		For curves on (very general) hypersurfaces of general type, what one has are a number
		of \emph{non-existence} results due to Voisin, Clemens and others (see \cite{clemens-ran}), 
		giving a lower bound on the genus in terms of the degree, but to my knowledge no
		existence results for good curves.\par
		A naive way to construct pairs $(C, X)$ is to start with a good curve $C\subset\P^n$ and let $X$
		be a hypersurface of sufficiently large degree containing $C$, which may be assumed smooth.
		Generally, the normal bundle $N_{C/X}$ has $\chi<0$, hence $H^1\neq 0$.
		The question then arises as to what good properties the pair $(C, X)$ has. This is addressed in
		some cases by our main result.

\subsection{New results: Stable balance, rigid regularity}
	Our aim here is to study normal bundles from the point of view of stability, regularity and rigidity.
	Motivated by the case of rational curves whose normal bundle, though balanced,
	 is not semistable unless it has integral slope, 
	we introduce a property weaker than stability 
	that we call \emph{stable balance}. This stipulates
	that the slopes of graded pieces of
	 the Harder-Narasimhan or stable filtration of $N$ should all lie
	in  some interval of length 1 (which must then contain the slope of $N$ itself). 
	Stable balance is obviously equivalent to balance for rational curves.
	For higher genus curves stable balance is weaker than interpolation.
Its advantage	is that it is meaningful for any bundle regardless of positivity or sections.
Because interpolation results are known for $\P^n$ and Fano hypersurfaces, we will
focus here on the case of Calabi-Yau and general type hypersurfaces, i.e. those of
degree $\geq n+1$ in $\P^n$. For curves on such hypersurfaces $X$, one usually has $\chi(N)<0$,
so 
the best can hope for as regards $N$  is that $C$ is rigid in $X$, i.e. $H^0(N)=0$.
And as for 'regularity', when $H^1(N)\neq 0$ the best kind of regularity one can hope for is 
what we call \emph{rigid-regularity} for the 
pair $(C, X)$ (see \S \ref{methods}), which is expressed in terms of vanishing for a suitable 'joint normal sheaf'
$N_{C/X/\P^n}$ in addition to rigidity of $C$ in $X$. 
Rigid-regularity implies that the pair is regular as such, i.e. deforms in a smooth family of the expected dimension,
and  this family projects injectively to the family of hypersurfaces.

\par
	Our main result is the following  (see \S\ref{results}):
	  \begin{thm}\label{main-thm}
	  	Assume
	  	\[ d\geq n+1\geq 5,\ \textrm{and}\]
	  	\[ \textrm{either\ } (g\geq 2 \textrm{\ and\ } e>(2n-4)(g-1)) \textrm{\ or\ } (g=1 \textrm {\ and\ } e>2n+2),
	  	\ \textrm{and}\]
	  	 \eqspl{main-display}{ 
	 \binom{d+n-2}{n-1}>ed/2.
	  } Then for a general  pair
	  	$(C,X)$, where $C$ is a general smooth curve of genus $g$ and degree $e$ and
	  	$X$ is a smooth  hypersurface of degree $d$  in $\P^n$,  
	  	$C$ is rigid-regular stably balanced on $X$.
	  	\end{thm}
	  	\textbf{Remarks:}\begin{enumerate}
	  	\item A result of Ballico \cite{ballico-twisted} allows one to improve the lower bound on $e$
	  	to  $e\geq (n-3)g+2n+6$ for $n\geq 8$.
	  	\item  	Note that  because $d>n$, we have $\binom{d+n-2}{n}\geq\binom{d+n-2}{n-1}$, 
	  	hence our assumptions imply
	  	\eqspl{binom}{\binom{d+n}{n}\geq\binom{d+n-1}{n}=\binom{d+n-2}{n-1}+\binom{d+n-2}{n}\\
	  	\geq ed/2+ed/2+1-g=ed+1-g.}
	  	Thus, the fact that $C$ has maximal rank implies that the restriction map
	  	$H^0(\O_{\P^n}(d))\to H^0(\O_C(d))$ is surjective.
	  \item
	  	Ideally,  one would like to replace the assumption \eqref{main-display}  by a weaker one such as
	  	$\binom{d+n}{n}>ed+1-g$, however our argument as is does not yield this. Likewise, our lower bound on
	  	$e$ might be weakened to something like nonspeciality or better yet, to the Brill-Noether number
	  	being nonnegative.
	\item
	  	Assumptions as in Theorem \ref{main-thm}, Case $g\geq 2$, in the family of pairs $(C, X)$,
	  	a general  $C$ has
	  	maximal rank and the family  clearly irreducible of the expected dimension; 
	  	hence so is the corresponding family of hypersurfaces
	  	$X$. As subfamily of the family of all hypersurfaces, the latter family has codimension
	  	\[h^1(N_{C/X})=-\chi(N)=(d-n-1)e+(n-4)(g-1).\]

	 
\item	
Here
			\[\chi(N)=e(n+1-d)+(n-4)(1-g).\]
	In the cases considered, this is always $\leq 0$.\item
Because $C$ is nonspecial
	and of maximal rank by \cite{ballico-ellia-max-rank}, the family of pairs $(C,X)$ is irreducible.\item
As stated above, earlier examples of rigid curves on hypersurfaces included a few Calabi-Yau
cases and general type examples satisfying $\chi(N)\geq 0$, 
which were not shown to be relatively regular 
or stably balanced. To my knowledge the above are the first examples with $\chi(N)<0$ 
( a property which makes deformation/smoothing arguments
trickier, and which leads us to consider deformations of the triple $C\subset X\subset \P^n$).\item
The range of hypersurface degrees for given $C$ is not necessarily optimal- an optimal one might be
$\binom{d+n}{n}>ed+1-g$.
Note that the degree, hence also the genus, of any rigid-regular curve of a hypersurface
 of degree $d$ is a priori bounded by a function of $d$. \par
\item
 It remains an open question in genus $g\geq 1$ whether the normal bundles
 to the curves we construct may actually be stable or semistable rather than
 just stably balanced. However Example \ref{unstable-ex} suggests that this may be
 false for stability ($n\geq 4$) or semi-stability ($n\geq 5$), at least for the degenerate curves that we construct. 
 By contrast, in the Fano case
 the normal bundle is often semistable (cf. \cite{semistable}).
\item
As for (in)completeness of our results, 
while there are numerous results from several viewpoints on the family of 
all curves on Fano hypersurfaces in $\P^n$, 
the situation changes drastically as soon as one enters the CYGT range $d> n$. 
In the very first such case, namely rational curves on the  general quintic in $\P^4$,
a well-known conjecture of Clemens stating that all such curves must be rigid (equivalently, regular, equivalently
have semistable normal bundle), has remained open after some 30 years despite many attempts.
Thus, anything close to completeness seems out of reach.
Our purpose here is merely to construct some good examples beyond the Fano range including some in
the general-type range, with no claim to sharpness or completeness.
\end{enumerate}
  
  \begin{example}
  	The Theorem inplies that there are
  	stably balanced,  rigid-regular curves of every degree $ 5<e<72$ and genus $0< g<e/5-5$
  	on  smooth quintic fourfolds.
  	\end{example}
\subsection{Methods}\label{methods}
We will use the fan  method developed in the author's earlier papers \cite{caudatenormal}. \cite{elliptic}.
After recalling some facts on fans, their hypersurfaces and degenerations
in \S \ref{fan-sec},
some properties of, and operations on stably balanced bundles are developed in \S \ref{stable-balance-sec}.
In particular, we study elementary modifications and specialization/generalization.
Next in \S\ref{higher-sec}, we develop some generalities on embedded deformations of a pair
$(C, X)$ in an ambient space such as $\P^n$, with focus on cases where $\chi(N_{C/X})<0$.
We define a suitable normal 
sheaf $N_{C/X/\P^n}$ controlling embedded deformations of the pair, such that
the vanishing of $H^1(N_{C/X/\P^n})$ is equivalent to the pair $(C,X)$ having unobstructed deformations
in $\P^n$ of the expected dimension $h^0(N_{C/X/\P^n})=\chi(N_{C/X/\P^n})$ 
(the cohomology of degree $>1$ usually vanishes automatically). The curve $C\subset X$ is
said to be \emph{rigid-regular} if $H^1(N_{C/X/\P^n})=0$ and $H^0(N_{C/X})=0$.
This implies that the space of deformations of the pair $(C,X)$ in $\P^n$, the so-called
flag Hilbert scheme, maps isomorphically to a smooth subvariety of codimension $h^1(N_{C/X})$ in
the deformation space of $X$. \par
The main result, Theorem \ref{main-thm}, presented in \S \ref{results},
has been described above.
   Its proof is essentially inductive in nature, based on a suitable fan degeneration similar to \cite{caudatenormal}, and
requires special care due to nonvanishing of  $H^1(N_{C/X})$. It relies in an essential way on
 maximal rank theorems of Hartshorne- Hirschowitz \cite{hh-droites} for general
 skew lines and
 Ballico-Ellia \cite{ballico-ellia-max-rank} for general nonspecial curves.
\par

 \subsection*{Conventions} We work over $\C$. As
 general references for deformation theory, see the books by Sernesi \cite{sernesi}
 or Hartshorne  \cite{hartshorne-def}.
 \par
 The slope of a bundle $E$ on a curve, i.e. $\deg(E)/\rk(E)$ is denoted by $\mu(E)$.
 \subsection*{Glossary of normal sheaves and corresponding deformations}
 
 \begin{itemize}\item
 $N_{C/X}$: deformations of $C$ in $X$.\item
 $T_X\mlog{C}$: deformations of pair $(C,X)$.\item
 $N_{C/X/P}$: deformations of pair $(C,X)$ in $P$.\item
 $N^C_{X/P}$: deformations of $X$ in $P$ fixing $C$.
 \end{itemize}
 \section{Fan and quasi-cone degenerations}\label{fan-sec}

 See \cite{caudatenormal} for details.
 We recall that a \emph{fan} (also called a 2-fan) is a reducible normal-crossing variety
 of the form
 \[P_0=P_1\cup_E P_2\]
 where
 \[P_1=B_p\P^n=\P_{\P^{n-1}}(\O(1)\oplus\O), P_2=\P^n\]
 and $E=\P(\O)\subset P_1$ is the exceptional divisor and $E\subset P_2$
 is a hyperplane. The family
 \[\mathcal P=B_{(p,0)}(\P^n\times\A^1)\]
 is called a standard fan degeneration and realizes $P_0$ as the
 special fibre in a family with general fibre $\P^n$.\par
 A \emph{hypersurface of type $(d_1, d_2)$} for $ d_1\geq d_2$ in $P_0$ has the form
 \[X_0=X_1\cup_Z X_2\]
 where \[X_1\in |d_1H_1-d_2E|_{P_1}, X_2\in |d_2H_2|_{P_2}, X_1\cap E=X_2\cap E=Z\]
 	($H_1, H_2$ are the respective hyperplanes). If $d_2=d_1-1$, $X_0$ is said to
 	be of \emph{quasi-cone type}. Given a family $\bar\cX\subset\P^n\times\A^1$
 	of hypersurfaces of degree $d_1$ whose special fibre has multiplicity $d_2$
 	at $p$, its birational transform $\cX\subset\cP$ is a family of
 	hypersurfaces in $\P^n$ specializing to one of type $(d_1, d_2)$. 
 	\section{Stably balanced bundles}\label{stable-balance-sec}
 	 \subsection{Basics}
 Given a vector bundle $E$ on a curve $C$, we denote by 
 \[(E_\bullet)=(0\neq E_1\subsetneq E_2...)\] 
 its Harder-Narasimmhan 
 (HN) filtration,
 characterized by the property that 
 \[\mu_{\max}(E)=\mu_1(E):=\mu(E_1)\] is the largest slope among 
 slopes of subbundles $E_1\subset E$,
 and $E_1$ is maximal of slope $\mu_1(E)$ (i.e. $E_1$ is the
 maximal-rank, maximal-slope subbundle and as such is unique and semi-stable); and likewise
 for $\mu_{i+1}(E):=\mu(E_{i+1}/E_i), \forall i\geq 1$. 
 Thus we have a strictly decreasing slope sequence
 \[\mu_{\max}(E)>...>\mu_{\min}(E)\]
( unless $E$ is semi-stable, in which case $\mu_{\max}(E)=\mu_{\min}(E)).$
 The interval $[\mu_{\min}(E), \mu_{\max}(E)]$ is called the slope interval
 of $E$ and its width $\mu_{\max}(E)-\mu_{\min}(E)$ is called the slope width.
 $E$ is said to be \emph{stably balanced} 
 if its slope width is $\leq 1$, i.e.
 \[\mu_{\max}(E)-\mu_{\min}(E)\leq 1.\]
 Note that the slope $\mu(E)$ is contained in the slope interval, hence if $E$
 is stably balanced then
 \[\mu_{\max}(E)-1\leq \mu(E)\leq \mu_{\min}(E)+1.\]
 By comparison, recall that $E$ is said to be \emph{balanced} if, denoting
 by $D_t$ a general effective divisor of degree $t$, we have for each $t\geq 0$ that either
 $H^1(E(-D_t))=0$ or $H^0(E(-D_t))=0$. Because this condition involves just negative
 'twists' of $E$, it makes sense just for 'positive' bundles, unlike the stable balance
 condition. It is also stronger than stable balance:\par
 	\begin{lem}\label{stable-balance-lem}
 	If $E$ is balanced then $E$ is stably balanced.
 	\end{lem}
 	\begin{proof}
 	Let $t$ be largest so that $H^1(E(-D_t))=0$ for $D_t$ general of degree $t$.\par
 	Then we have $H^1(E'(-D_t))=0$ for any quotient $E'$ of $E$.
 	Also  $H^0(E(-D_{t+1}))=0$ and    $H^0(E"(-D_{t+1}))=0$ for any subsheaf $E"\subset E$.
 	Therefore $\chi(E(-D_t))\geq 0$ and ditto for a minimal slope quotient,
 	 while $\chi((E(-D_{t+1})))\leq 0$ and ditto for the max slope  subsheaf.
 	This easily implies $\mu_{\max}(E)+1-g\leq t+1$ and similarly $\mu_{\min}(E)+1-g\geq t$,
 	so $\mu_{\max}(E)-\mu_{\min}(E)\leq 1$.
 	\end{proof}
 	In fact the same argument yields the following stronger
 	 result, noted by Vogt \cite{vogt-interpolation}:
 	\begin{lem}
 	If $E$ is balanced then we have for every subbundle $F\subset E$,
 	\[\mu(F)\leq\ceil{\mu(E)}\]
 	In particular, if $E$ is balanced and  $\mu(E)$ is an integer then $E$ is semistable. 
 	\end{lem}
 	\begin{proof}
 	Adding $1-g$ to both sides, it suffices to prove
 		$\chi(F)/\rk(F)\leq\ceil{\chi(E)/\rk(E)}$.
 	If $t\geq\chi(E)/\rk(E)$ then $\chi(E(-D_t))\leq 0$, hence $H^0(E(-D_t))=0$,
 	hence $H^0(F(-D_t))=0$, hence $\chi(F(-D_t))\leq 0$, i.e.
 	$t\geq\chi(F)/\rk(F)$.
 	\end{proof}
 	Lemma \ref{stable-balance-lem}  implies that the (high-degree) 
 	balanced curves constructed, e.g. in \cite{alyang}, \cite{caudatenormal}
 	and \cite{elliptic}, all with positive normal bundle in Fano hypersurfaces, are all stably balanced.
 	So the interest in stably balanced curves on hypersurfaces is mainly for Calabi-Yau
 	and general-type hypersurfaces (degree $d\geq n+1$ in $\P^n$, 
 	as well as in lower-degree curves on Fano hypersurfaces ($d\leq n$) .
 
A semistable bundle of slope $<0$ has $H^0=0$.  
 By Serre duality, a semistable bundle of slope $\mu>2g-2$ has $H^1=0$.
 Applying this to the gradeds of the HN filtration, we conclude
 \begin{lem}\label{vanishing-lem}
 Suppose $E$ is stably balanced of slope $\mu$. If
 $\mu<-1$ then $H^0(E)=0$. If $\mu>2g-1$ then $H^1(E)=0$.
 \end{lem}
 \begin{example}
If $C$ is a $(g,e)$ curve in $\P^n$ with stably balanced normal bundle $N$, such that
$e>1+(n-2)g$, then $H^1(N(-1))=0$, which means geometrically that for a 
hyperplane $H$, and $C$ general in its (smooth) Hilbert scheme component,
$H.C$ consists of $e$ general points on $H$.\par
  If $C\simeq \P^1$, stable balance is equivalent to balance.\par
\end{example}
 It is noted in \cite{elliptic}, Lemma 8 that if the gradeds of the HN filtration of
 $E$ are direct sums of line bundles with degrees in some interval
 $[a, a+1], a\geq 2g-1$ (so that $E$ is stably balanced), then $E$ is balanced.
 Otherwise, the relation between stable balance and balance in the sense
 used, e.g. in \cite{elliptic} is not understood for $g>0$.
 	\par
 	\subsection{Modifications}
 Our first result is about the behavior 
 of stable balance under general modifications. This was essentially
 done in \cite{caudatenormal} for $C=\P^1$ and
  the discussion there largely extends to the general case.
  This result will be our main tool in proving stable balance for the normal bundle of the curves
  that we construct for the Main Theorem.
  \par
 By a \emph{semi-general} (resp. \emph{general})
  down modification of $E$ we mean an exact sequence
 \[\exseq{E'}{E}{\tau}\]
 such that $\tau$ is of the form $\bigoplus\kk(p_i)^{s_i}$
 for some collection (resp. a general collection)
 of distinct points $p_1,...,p_k\in C$, and such that
 the map
 \[E|_{p_i}\to\kk(p_i)^{s_i}\]
 is a general quotient for $i=1,...,k$. 
 \begin{lem}\label{mod-lem}
 	A semi-general down modification of a stably balanced bundle is stably balanced.
 \end{lem}
 \begin{proof}
 	We may assume the modification $E\to \tau$ is
 	at a single  point $p$.\par
 	Note that $\mu_{\max}(E')\leq\mu_{\max}(E) $.\par
 	Now assume to begin with that $E$ is semistable. Then for any quotient
 	$E'/G$, we have 
 	\[\mu(E)\leq\mu(E/G)\leq\mu(E'/G)+1.\]
 	Thus $\mu_{\min}(E')\geq\mu(E)-1$, so $E'$ is stably balanced.\par
 	In the general case consider the HN filtration
 	\[0=E_0\subsetneq E_1\subsetneq...\subsetneq E_a=E.\]
 	and the corresponding filtration $E'_\bullet$ on $E'$. There is a unique $k$ such that
 	the map $E_i|_p\to \tau$ is injective for $i\leq k$ and surjective, non-injective for $i>k$. The
 	we have
 
 	\[E'_i=E_i(-p), \mu(E'_i/E'_{i-1})=\mu(E_i/E_{i-1})-1, i<k,\]
 	\[\mu(E_k/E_{k-1})-1\leq\mu_{\min}(E'_k/E'_{k-1})\leq\mu_{\max}(E'_k/E'_{k-1})\leq\mu(E_k/E_{k-1})\]
 	\[\mu(E'_i/E'_{i-1})=\mu(E_i/E_{i-1}), i>k.\]
 	Now $E'$ admits a filtration, not necessarily the HN filtration, with semistable gradeds consisting of
 	\[E_i(-p)/E_{i-1}(-p), i=1,...,k-1;  E_i/E_{i-1}, i=k+1,...,a;\] 
 	plus the HN graded pieces of $E'_k/E'_{k-1}$.
 	By the semistable case just discussed, 
 	the slopes of these all lie in the interval $[\mu(E_k/E_{k-1})-1, \mu(E_k/E_{k-1})]$,
 	hence $E'$ is stably balanced. 
 	\end{proof}
 	Thus, semi-general elementary modifications stay within the class of stably balanced 
 	bundles. As the following results shows, this is false for the class of (semi)stable bundles.
 	\begin{example}\label{unstable-ex}
 	In any genus, if $E=2\O$ and $\tau=\kk(p)$, then $E$ is semistable but $E'=\O\oplus\O(-p)$
 	is stably balanced but not semistable. So in that sense the result is sharp.\par
 	Here we give examples in any genus $g\geq 2$ where $E$ is stable but its modification $E'$ is not.
 We	start with an example of a non-stable but maybe semi-stable
 modification of a stable rank-3 bundle. Consider a general extension
 	\[0\to{\O}\stackrel{s}{\to}{F}\to{G}\to 0\]
 	where $G$ is a  line bundle of degree $1$ 
 	and the extension corresponds to a general element  $\epsilon\in H^1(  G\inv)=H^0(G\otimes K)^*$. 
 	Here $G\inv$ denotes the dual or inverse of the line bundle.
 	I claim that $F$ is stable. In fact, consider any line bundle $L$ of degree $\geq\mu(F)=1/2$
 	i.e. $\deg(L)\geq 1$.
 	. If $\deg(L)>1$
 	or $\deg(L)=1$ and $L\neq G$ then clearly $H^0(F\otimes L\inv)=0$. If $L=G$, the same
 	is true by nontriviality of the extension. 
 	Thus $F$ is stable.
\par 	
Now consider an extension
 	\eqspl{E}{\exseq{\O}{E}{F}}
 	corresponding to an  element $\eta\in H^1(\check F)=H^0(F\otimes K)^*$
 	which comes from $H^1( G\inv)$, i.e. fits in an exact diagram
 	\eqspl{}{
 	\begin{matrix}
 	0\to&\O&\to&E&\to&F&\to0\\
 	&\downarrow&&\downarrow&&\downarrow&\\
 	0\to&\O&\to&E_G&\to&G&\to 0.
 	\end{matrix}
 	}
 	Let $L$ be any line bundle with $\deg(L)\geq\mu(E)=1/3 $
 	(i.e. $\deg(L)>0$). As we have seen $H^0(F\otimes L\inv)=0$,
hence $H^0(E\otimes L\inv)=0$. So $E$ has no rank-1 destabilizing subsheaf.
 \par
 Now we have an exact sequence
 \[\exseq{F}{\wedge^2E}{G}\]
 which is none other than the dual of \eqref{E}, twisted by $G$,
and we know $H^0(F\otimes L\inv)=0$. Therefore 
an argument similar to the above
shows that 
we have  $H^0(\wedge^2E\otimes L\inv)=0$, for any  $L$ of degree $>0$,
so $E$ has no rank-2 destabilizing subsheaf.
Thus $E$ is stable (even cohomologically stable in the sense of
\cite{ein-laz-normal} ).\par
 	Now a 
 	general corank-1 modification $E'$ of $E$
 	at $p$ fits on an exact sequence
 	\[\exseq{\O(-p)}{E'}{F}\]
 	which corresponds to an element of $H^1(\check F\otimes\O(-p))$.
 As with the original element $\eta$ this comes from an element of
 	$H^1(\check G\otimes \O(-p))$ 	corresponding to a modification $E'_G$ of $E_G$.
 	Therefore the section $s$ of $F$ lifts to $E'$.
 	Then $E'$ has slope $0$ but $H^0(E')\neq 0$ so $E'$ is not stable
 	(but maybe semi-stable).\par
 	For an example of a non-semistable modification assume $g\geq 3$ and
 	consider instead a similar extension
 	\eqspl{E2}{\exseq{2\O}{E_2}{F}
 	}
 	which may also be realized as an extension
 	\[\exseq{3\O}{E_2}{G}\]
 	corresponding to  a general element of $H^1(3G\inv)$.
 	then a similar argument as above shows that for any line bundle $L$ with
 	$\deg(L)\geq\mu(E_2)=1/4$, i.e. $\deg(L)\geq 1$,
 	we have $H^0(E_2\otimes L\inv)=0$. Also, 
 	\[H^0(\wedge ^3E_2\otimes L\inv)=H^0(\check E_2\otimes G\otimes L\inv)\]
 	which clearly vanishes when $\deg(L)>0, L\neq G$, and also for $L=G$ due to injectivity
 	of $H^0(3\O)\to H^1(G)\simeq H^1(\O)$ by generality of the extension.
 	Moreover we have an exact sequence
 	\[\exseq{\O}{\wedge^2E_2}{E_1}\]
 	where $E_1$
 is analogous to $E$ above. Then a similar argument again shows $H^0(\wedge^2E_2\otimes L\inv)=0$.
 Therefore $E_2$ cannot admit a destabilizing subsheaf of rank 1,2 or 3,
 hence it is stable.\par
 Now a general corank-2 modification at $p$, $E_2'$, fits in an exact sequence
 \[\exseq{2\O(-p)}{E'_2}{F}\]
 and as above we have $H^0(E'_2)\neq 0$ while $\mu(E'_2)=-1/4$, so $E'_2$
 is not semistable.\qed
 	\end{example}
 	
%
%
%

 \subsection{Degeneration}
Next we come to an important technique for constructing stably balanced bundles
via degeneration. We begin with the following situation:
\[\pi:\mathcal C\to B\]
is a family of curves over a smooth curve with all fibres $C_t$ smooth except
\[C_0=C_1\cup_pC_2\]
with $C_1, C_2$ smooth and transverse; $\mathcal E$ is a vector bundle on
$\mathcal C$, and we set $E_i=\mathcal E|_{C_i}, i=1,2/$, $E_t=E_{C_t}$.
\begin{lem}
	Notations as above, assume $E_1$
	is semistable and $E_2$ is stably balanced. Then
	a general $E_t$ is stably balanced.
	\end{lem} 
\begin{proof}
	Let $[\alpha, \beta]$ be the slope interval of $E_2$ and let $F_t$
	be a subbundle of the general fibre $E_t$ ( $F_t$ may possibly be
	 defined only after
	a base-change, though actually we may choose for $F_t$ the maximal-rank,
	maximal-slope subbundle and then no base-change is required). 
	We will show that \eqspl{sub-eq}{\mu(F_t)\leq \mu(E_1)+\beta.}
	Then applying this result to the dual will show that
	\[\mu(F_t)\geq\mu(E_1)+\alpha.\] Thus,
	the slope interval of $E_t$ is contained in $[\mu(E_1)+\alpha, \mu(E_1)+\beta]$,
	proving our result.\par
	Now to prove \eqref{sub-eq},  after a suitable base-change $B'\to B$ and replacing $\mathcal C$
	by $\mathcal C'=\mathcal C\times_BB'$, we may assume given a torsion-free subsheaf
	$\mathcal F\subset\mathcal E$ extending $F_t$. Then after further base-change and
	blowing up we may assume 
	 $\det(F_t)=\wedge^sF_t$ extends to an invertible
	sheaf $\mathcal L$ on $\mathcal C'$ and that the map 
	$\mathcal L\to\wedge^s\mathcal E$ vanishes only on an effective divisor
	$Z$ supported on the special fibre so the induced map $\mathcal L(Z)\to
	\wedge^s\cE$ is nowhere vanishing. Then we have
	\[s\mu(F_t)=\deg(\mathcal L(-Z)|_{C'_0})\leq\deg(\mathcal L_{C_1})
	+\deg(\mathcal L|_{C_2})\leq s(\mu(E_1)+\beta).\]
	
\end{proof}
The same result with the same proof also holds for a comb-shaped curve:
\begin{lem}\label{comb-lem}
	Let $C_0$ be a connected nodal curve of the form 
	$D\cup C_{01}\cup...\cup C_{0m}$
	with $D\cap C_{0i}$ a single point, $i=1,..,m$, and no other singularities.
	Let $\cC\to B$ be a family as above with special fibre $C_0$ and
	general fibre $C_t$, and let  $\cE$
	a vector bundle on $\cC$ such that 
	the restriction of $\cE$ on all but one of the components $D, C_{01},...,C_{0m}$
	is semistable, and the restriction is stably balanced on the remaining component.
	Then the restriction $\cE|_{C_t}$ is stably balanced. 
	\end{lem}
\begin{rem}	Note that a stably balanced bundle $E$ of slope $\mu(E)\leq -1$
	must have $\mu_{\max}(E)<0$, hence $H^0(E)=0$.\par
	A lci subvariety $C\subset X$ is said to be stably balanced if the normal bundle
	$N=N_{C/X}$ is, and rigid if $H^0(N)=0$. When $C$ is a curve of genus $g$ 
	and $X$ is $n$-dimensional, we have
	\[\mu(N)=(-C.K_X+2g-2)/(n-1).\]
	Consequently, if $C$ is stably balanced and $2g-2\leq C.K_X-n+1$ then $C$ is rigid.
\end{rem}\subsection{Relation to interpolation}	
	Lemma \ref{vanishing-lem} implies that a general divisor
	$D_t=p_1+....,+p_t$ on $ C$,  we have $H^1(N(-D_t)=0$ provided
	$t<\mu(N)-(2g-1)$. Therefore we conclude:
	\begin{lem}\label{interpolation-lem}
	If $C\subset X$ is a stably balanced curve of genus $g$ and $t<(-C.K_X+2g-2)/(n-1)-2g+1$,
	the $C$ is $t$-interpolating in the sense that $t$ general points on a general deformation
	of $C$ are general in $X$.
	\end{lem}	The $t$-interpolating property in the Lemma is weaker that the full interpolation or balance property as
		studied e.g. in \cite{elliptic}.\par
\section{Generalities on deformations}\label{higher-sec}
Here we consider curves $C$ on hypersurfaces $X$ of degree $d\geq n+1$ in $\P^n$.
This case is of a different character  compared to the case $d\leq n$
because generally $\chi(N_{C/X})<0$ so
$H^1(N_{C/X})$ cannot vanish, nor can $X$ be general. Consequently
at least the smoothing issue has to be handled differently. We shall
do so by studying deformations and specializations of the triple $(C,X,\P^n)$.
\subsection{Deformations of pairs} Compare  \cite{sernesi}.
Let $C\subset X$ be an inclusion of smooth projective varieties. There are corresponding deformation
functors and spaces
\[\Def(C/X)\to \Def(C, X)\to\Def(X)\]
parametrizing deformations of, respectively, $C$ fixing $X$, the pair $(C, X)$, the variety $X$.
These correspond to an exact sequence of sheaves on $X$
\[\exseq{T_X\mlog{C}}{T_X}{N_{C/X}}.\]
In fact, $T_X$ is a sheaf of Lie algebras and $T_X\mlog{C}$ is a subalgebra sheaf, which endows 
$N=N_{C/X}$ with the structure of Lie atom (see \cite{atom}), which yields the corresponding deformation theory,
i.e. deformations of $C$ in a fixed $X$. Thus the complex $N[-1]$ is endowed with a bracket
in the derived category and deformations (resp. obstructions  to deformations)
 of $C$ fixing $X$ are in $H^0(N)$ (resp. $H^1(N)$); for the pair $(C, X)$ the corresponding groups
 are $H^1(T_X\mlog{C}), H^2(T_X\mlog{C})$.
 
 \subsection{Embedded deformations of pairs}\label{embedded-sec}
 Let $X\subset P$ be a locally complete intersection with $P$ smooth, and let $C\subset X$
 be a locally complete intersection in $X$.
 Then the normal bundles $N_{C/P}$ and $N_{X/P}$ admit $N_{X/P}|_C$ as a common quotient. 
 Define a coherent sheaf $N_{C/X/P}$ on $X$ by the exact sequence
 \eqspl{ncxp}{\exseq{N_{C/X/P}}{N_{C/P}\oplus N_{X/P}}{N_{X/P}|_C}.}
 Then $N_{C/X/P}$ controls deformations of the pair $(C, X)$ as subvarieties of $P$,
 and we have an exact sequence
 \[\exseq{N_{C/X}}{N_{C/X/P}}{N_{X/P}}.\] 
 We will say that the pair $(C\subset X)$ of subvarieties is 
 \emph{rigid-regular} in $P$ provided the
following  vanishings hold.
 \[	H^0(N_{C/X})=0, H^1(N_{C/P/X})=0, H^1(N_{X/P})=0.\]
The following follows from general deformation/obstruction theory (see e.g. \cite{sernesi}):
 \begin{lem}
 Assume $(C,X)$ is rigid-regular in $P$.
 	Then near $(C,X)$, the Hilbert scheme of pairs in $P$ embeds as a smooth subvariety
 	of codimension $-\chi(N_{C/X})=h^1(N_{C/X})$ in the Hilbert scheme of $X$ in $P$.
 	\end{lem}
 There is another exact sequence
 \eqspl{normal-seq}{\exseq{N^C_{X/P}}{N_{C/X/P}}{N_{C/P}}}
 where $N^C_{X/P}\subset N_{X/P}$ is the subsheaf corresponding to deformations
 of $X$ fixing $C$ (or equivalently, deformations of the birational
 transform of $X$ in the blowup of $P$ in $C$).\par 
 Now suppose $X$ is a hypersurface of degree $d$ with equation $f$ in $P=\P^n$. 
 We can identify
 \eqspl{normal-ident}{N_{X/P}=\O_X(d), N_{X/P}^C=\I_{C/X}(d)=\I_{C/P}(d)/f\O_P.}
 \par Now assume $C, X, P$ are smooth, and let $T_{P/X/C}\subset T_P$ denote the
 subsheaf consisting of local vector fields tangent to both $X$ and $C$.
 Then we have an exact sequence
 \[\exseq{T_{P/X/C}}{T_P}{N_{C/X/P}}\]
 and also
 \[\exseq{T_P(-X)}{T_{P/X/C}}{T_{X}\mlog{C}}.\]
 Then from the cohomology sequence we conclude:
 \begin{lem}\label{regularity-lem}
 (i)	Assume the vanishings
 	\[H^1(N_{C/X/P})=0, H^2(T_P)=0, H^3(T_P(-X))=0.\]
 	Then $C$ is rigid-regular in $X$, i.e.  
 	the pair $(C,X)$ has unobstructed deformations of dimension $h^1(T_{X}\mlog{C})$.\par
 (ii)	If moreover $H^1(N^C_{X/P})=0$, then $C$ may be assumed general in the Hilbert scheme 
 of subvarieties $C\subset P$.
 	\end{lem}
 	\begin{rem}\label{max-rk-rem}
 	If $C$ is a smooth curve and $P=\P^n$, the vanishing $H^1(N^C_{X/P})=0$ follows from 
 	maximal rank, specifically from  $H^1(\I_{C/P}(d))=0$.
 	\end{rem}
 \begin{example}
 	Let $P=\P^n, n\geq 5$ and $X\subset P$ a smooth hypersurface
 	of degree $d$. Then if $C\subset X$ is a curve of degree $e$ and genus $g$
 	with  $H^1(N_{C/X/P})=0$,
 	the pair $(C, X)$ has unobstructed deformations of dimension $h^1(T_X\mlog{C})$.
 	If $(C,X)$ is relatively regular in $P$, the family of hypersurfaces carrying 
 	deformations of $C$ is smooth of codimension $e(d-n-1)+(n-4)(g-1)$ near $X$,
 	and surjects to the Hilbert scheme of deformations on $C$ in $P$, and in particular $C$
 	has maximal rank or more specifically $h^1(\I_{C/P}(d))=0$.
 	\end{example}
%
 \subsection{Singular case}\label{singular-sec}
 Consider a lci diagram of connected normal-crossing varieties
 \[C_0=C_1\cup C_2\subset X_0=X_1\cup_Z X_2\subset P_0=P_1\cup_Q P_2\subset\mathcal P\]
 where each $C_i, X_i, P_i$ and $\mathcal P$ are smooth and there is a morphism
 $\pi:\mathcal P\to B$ to a smooth curve such that $P_0=\pi\inv(0)$.
 Then a general fibre $P_t=\pi\inv(t)$ is smooth.  We have an exact sequence
 \[\exseq{N_{C_0/X_0/P_0}}{N_{C_0/X_0/\mathcal P}}{\O_{X_0}}.\]
We have an isomorphism $H^0(\O_{X_0})\simeq H^0(\O_Z)$
and $\O_Z$ can be identified with the  $\mathcal T^1_{P_0}\otimes\O_{X_0}=\mathcal T^1_{X_0}$ which
 restricts on $C_0$ to $\mathcal T^1_{C_0}$, it follows that a deformation
 of $(C_0, X_0)$ in $\mathcal P$ is a smoothing of $C_0$ and $X_0$ iff its image in $\O_{X_0}$
 is nonzero. Then from the long cohomology  sequence we conclude
 \begin{lem}
 Notations as above, assume 
 \[H^1(N_{C_0/X_0/P_0})=0=H^1(\O_{X_0}).\] Then
 $(C_0, X_0, P_0)$ is smoothable to $(C_t, X_t, P_t)$ in a general fibre $P_t$ of $\pi$,
 i.e. there is a smooth variety $T$ of dimension $h^0(N_{C_0/X_0/\mathcal P})$,
 a smooth morphism $T\to B$ with a point $0'$ over $0$,
  and a triple $(\mathcal C, \mathcal X, \mathcal P\times_{B}T)$
 flat over $T$, with special fibre $(C_0, X_0, P_0)$ and general fibre $(C_t, X_t, P_t)$
 with $C_t, X_t, P_t$ smooth and moreover $H^1(N_{C_t/X_t/P_t})=0$.	
\end{lem}
As for the vanishing in question, note that 
\[N_{C_0/X_0/P_0}|_{X_i}=N_{C_i/X_i/P_i}, i=1,2\]
hence we have an exact sequence
\[\exseq{N_{C_1/X_1/P_1}(-Z)}{N_{C_0/X_0/P_0}}{N_{C_2/X_2/P_2}}.\]
Note finally that the isomorphism  $\mathcal T^1_{P_0}\otimes\O_{X_0}=\mathcal T^1_{X_0}$
implies that the smoothing $\cP$ of $P_0$ to which $X_0$ extends
is also a smoothing of $X_0$, hence
\begin{lem}\label{smoothing-lem} 
Notations as above, the general  fibre of $X/B$ is smooth.
\end{lem}

%
 \section{Proof of main results}\label{results}
 We will now prove Theorem \ref{main-thm} as stated in the Introduction..
%

\begin{proof}[Proof of Theorem] We use
 induction on $d\geq n+1$ and a suitable quasi-cone degeneration in a  fan degeneration of $\P^n$
 as in \S\ref{fan-sec} .
 We construct a suitable pair $(C_0, X_0)$ of a curve on a 
 on a quasi-cone degeneration as in \S \ref{singular-sec}
 \[C_0=C_1\cup C_2\subset X_0=X_1\cup_Z X_2\subset P_0=P_1\cup_QP_2.\]
 Thus $P_1=B_p\P^n=\P_{\P^{n-1}}(\O(1)\oplus\O)$ with exceptional divisor 
 $Q=\P_{\P^{n-1}}(\O)\simeq\P^{n-1}$, $P_2=\P^n$ with $Q\subset P_2$ a hyperplane,
$X_1\subset P_1$ is the blowup of a hypersurface of degree $d$
with multiplicity $d-1$ at $p$,  and $X_2\subset P_2$ is a hypersurface of degree $d-1$.
As noted above (via projection from $p$), $X_1$ may be constructed 
as $B_Y\P^{n-1}$ where $Y=F_{d-1}\cap F_d$
 is a $(d-1, d)$ complete
intersection to be specified below.
Note that if we show that the triple $(C_0, X_0, P_0)$ extends to
$(\cC, \cX, \cP)$ for $\cP$ a smoothing of $P_0$ then in the general fibre
$(C, X, P^n)$, $C$ and $X$ are smooth.
\par
We now assume $d\geq n+1$ and work by induction on $d$.
Assume to begin with that $e$ is even.
We first construct $C_2\subset X_2\subset P_2$. Begin with a general curve $C_2\in\P^n=P_2$
of genus $g\geq 0$ and degree $e_2=e/2>g+n$ satisfying the inequalities \eqref{main-display}.
Note that because $d\geq n+1$ we have $\binom{d+n-2}{n-1}\leq\binom{d+n-2}{n}$ hence

\eqspl{c2}{
	\binom{d+n-2}{n}>e_2(d-2)+1-g.
} 
By the Maximal Rank Theorem of Ballico-Ellia 
\cite{ballico-ellia-max-rank},
the restriction map
\[H^0(\O_{\P^n}(d-2))\to H^0(\O_{C_2}(d-2)h)\]
is surjective, i.e. $H^1(I_{C_2}(d-2))=0$. 
Now I claim that a general hypersurface of degree $d-1$ containing $C_2$ is smooth.
To that end, note \[H^2(I_{C_2}(d-3))=H^1(\O_{C_2}((d-3)h)=0,\]
so $I_{C_2}$ is $(d-1)$-regular, therefore $I_{C_2}(d-1)$ is globally generated. 
Consequently by Bertini's Theorem a general hypersurface $X_2\subset\P^n$ of degree $d-1$
containing $C_2$   will be smooth off $C_2$. As for smoothness on $C_2$, 
note that $H^0(I_{C_2}(d-1))$
generates $I_{C_2}(d-1)$ hence  generates 
its quotient $\check N_{C_2/\P^n}(d-1)$; and as the latter bundle has rank $n-1>1$, 
the image of a general element of  $H^0(I_{C_2}(d-1))$ will be nowhere-vanishing. 
This means thaat  a general hypersurface $X_2$ of degree $d-1$ containing $C_2$ will be smooth along
$C_2$ as well.\par
Now by induction if $d>n+1$ or, in the initial case $d=n+1$, by the results of
\cite{alyang} or  \cite[Corollary 28]{elliptic},  
the normal bundle $N_{C_2/X_2}$ is stably balanced (or, in case $d=n+1$, actually balanced,
hence by Lemma
\ref{stable-balance-lem}, stably balanced).
Then by Lemma \ref{vanishing-lem} and our assumed lower bound on $e$,
it follows that 
\[H^1(N_{C_2/P_2}(-h))=0.\] 
Alternatively for $n\geq 8$, 
this also follows, even with a smaller lower bound
on $e$, from Ballico's result \cite{ballico-twisted}.  From this, together with the basic exact sequence
\eqref{ncxp} it follows that
\[H^1(N_{C_2/X_2/P_2}(-h))=0\]
and also that by varying $C_2, X_2$ while fixing $Z=X_2\cap Q$, the points 
$W\setminus C_2\cap Z=\{p_1,...,p_{e_2}\}$
will be $e_2$ general points on $Z$. \par
Now to construct $C_1\subset X_1\subset P_1$, identify $Z$ with the
hypersurface $F_{d-1}\subset\P^{n-1}$ (of degree $d-1$),
and let $\l_i\subset\P^{n-1}$ be a general line through $p_i$. As $\{p_1,...,p_{e_2}\}$
is a general collection of points on $Z$ and each $\l_i$ is a general line through
$p_i$, $\{\l_1,...,\l_{e_2}\}$ is a general collection of lines in $\P^{n-1}$.
By assumption, we have \[\binom{d+n-2}{n-1}>e_2d.\]
Now by the Maximal Rank Theorem of Hartshorne and Hirschowitz \cite{hh-droites},
the restriction map
\[H^0(\O_{\P^{n-1}}(d-1))\to \bigoplus H^0(\O_{ L_i}(d-1))\]
is surjective.
Let $F_d\subset\P^{n-1}$ be a general hypersurface of degree $d$ containing
the $d-2$ points $\l_i\cap F_{d-1}\setminus p_i$ but not $p_i$ for $i=1,...,e_2$.
Now I claim that 
$F_d$ has general tangent hyperplanes at $F_{d-1} \cap \l_i\setminus p_i, i=1,...,e_2$. Indeed this is clear because
we may take for $F_d$ a general union of  $d$ hyperplanes where the first $d-2$
go through the $d-2$ points of $I_i\setminus p$.\par

Now let $Y=F_{d-1}\cap F_d$ and $X_1=B_Y\P^{n-1}\subset P_1$ as above.
Note that the lines $\l_i$ lift to conics $D_i$ on $X_1$ (i.e. birational transforms of conics  from $\P^n$)
and we set
\[C_1=D_1\cup...\cup D_{e_2}.\]To avoid confusion we point out that The Hartshorne-Hirschowitz
result was applied wo the image of $C_1$ in $\P^{n-1}$, not
to its image in $\P^n$.
Note that by the claim above, there is for each $i$
an open set $U_i$ of hypersurfaces $F_d$ such that the normal bundle $N_{D_i/X_1}$ 
is a general modification
of $N_{L_i/\P^{n-1}}$, hence it is balanced. 
Therefore, there is an open set $\bigcap U_i$  of hypersurfaces $F_d$ such that  $N_{D_i/X_1}$
is balanced for all $i$. In particular, then 
$H^0(N_{C_1/X_1}(-Z))=0$.

Now on $P_1$ we have, where $H$ is a hyperplane from $\P^n$ and $L_1=H-Q$ is a hyperplane from $\P^{n-1}$:,
\[X_1\in|dH-(d-1)Q|=|dL_1+Q|\]
hence \[N_{X_1/P_1}=\O(dL_1+Z)=\O((2d-1)L_1-E_Y).\]
Now note that $P_1=\P_{\P^{n-1}}(\O(L_1)\oplus\O)$ hence
$\pi_*(\O(Q))=\O\oplus\O(-L_1)$, $R^1\pi_*(\O(Q))=0$, hence
\[H^1(\O_{P_1}(dL_1+Q))=H^2(\O_{P_1})=0.\]
From the normal sequence
\[\exseq{\O_{P_1}}{\O_{P_1}(dL_1+Q)}{N_{X_1/P_1}}\]
 it then  follows easily that $H^1(N_{X_1/P_1})=0$. Also clearly $H^1(N_{C_1/ P_1})=0$. 
 Now consider the exact diagram
 \eqspl{}{
 \begin{matrix}
 0\to&\O(dL_1)&\to&\O(dL_1+Q)&\to&\O_Q((d-1)L_1)&\to 0\\
 &\downarrow&&\downarrow&&\downarrow&\\
 0\to&\O_{C_1}(dL_1)&\to&\O_{C_1}(dL_1+Q)&\to&\O_{C_1.Q}((d-1)L_1)&\to 0.
 \end{matrix}
 }

By the aforementioned Maximal Rank Theorem of Hartshorne and Hirschowitz \cite{hh-droites}, 
the left map induces a surjection on $H^0$ and clearly so is the case for the right map,
by the generality of the points $C_1.Q$ and  the bound $ e_2 \leq \binom{n-d+1}{d-1}$.
Therefore,  ditto for the middle map, i.e.
the restriction map
\[H^0(P_1, dL_1+Q)\to H^0(C_1, (dL_1+Q)|_{C_1})\]
is surjective. Hence
we conclude that
\[H^1(N_{C_1/X_1/P_1})=0.\]
Therefore by the obvious exact sequence
\[\exseq{N_{C_2/X_2/P_2}(-1)}{N_{C_0/X_0/P_0}}{N_{C_1/X_1/P_1}}\]
(or just as well, the analogous one with $C_1,C_2$ etc. interchanged), we conclude
\[H^1(N_{C_0/X_0/P_0})=0.\]
Next we claim that $C_0$ is rigid on $X_0$, hence likewise for the smoothing.
To see this note that
\[N_{C_0/X_0}|_{C_i}=N_{C_i/X_i}, i=1,2\]
and as we have seen $N_{C_1/X_1}|_{D_i}$ is a general modification of $N_{\l_i/\P^{n-1}}$
	and in particular $N_{C_i/X_1}|_{D_i}(-Z)$ is strictly negative 
	and its upper subspace at $p_i$ is general.
	As for $N_{C_2/ X_2}$, in the initial case $d=n+2$, we know it is balanced and we have
	$\chi(C_2/X_2)=(n-4)(1-g)<e_2$ so $h^0(N_{C_2/X_2})<e_2$. In case $d>n+2$ we have inductively
	that $h^0(N_{C_2/X_2})=0$. Thus in all cases it follows 
	from the obvious exact sequence
	\[\exseq{N_{C_1/X_1}(-Z)}{N_{C_0/X_0}}{N_{C_2/X_2}}\]
	that
\[H^0(N_{C_0/X_0})=0.\]
Moreover an argument as in the proof of \cite{caudatenormal}, Theorem 6
or \cite{cyrigid}, Theorem 1 (essentially, blowing down the tails
$D_i$)   shows that for a smoothing $C'\subset X'$ of $C_0\subset X_0$,
the normal bundle $N_{C'/X'}$ is essentially a deformation of a semi-general
down modification of $N_{C_2/X_2}$ corresponding to the upper subspaces
of $T_{p_i}D_i$ and consequently by Lemma \ref{mod-lem}, $C'\subset X'$ is stably balanced.\par
To complete the proof it suffices by Lemma \ref{regularity-lem} to prove that $H^1(N^{C'}_{X'/\P^n})=0$.
By Remark \ref{max-rk-rem}, it suffices to prove  that $H^1(\I_{C'/\P^n}(d))=0$.
For this it would suffice to prove that $H^1(\I_{C_0/P_0}(L_0))=0$ where $L_0$
is the following lne bunde on $P_0$:
\[L_0= \O_{P_1}(dH-(d-1)Q)\cup \O_{P_2}((d-1)H).\]
Using $H^1(\O_{P_2}((d-2)H))=0=H^1(\O_{C_2}((d-2))$ (the latter due to 
$e > 2g + 2n$), 
we have an exact diagram
\eqspl{restriction}{
\begin{matrix}
0\to&H^0(\O_{P_2}((d-2)H))&\to&H^0(L_0)&\to&H^0(\O_{P_1}(dH-(d-1)Q))&\to  0\\
&\downarrow&&\downarrow&&\downarrow&\\
0\to& H^0(\O_{C_2}(d-2)) &\to&H^0(\O_{C_0}(L_0))&\to&H^0(\O_{C_1}(dH-(d-1)Q) )& \to 0.
\end{matrix}
}
As we have seen, the left vertical map is surjective. As for the right vertical map, we can
argue as above:  identify
$\O(dH-(d-1)Q)=dL_1+Q$ and use the exact sequence
\[\exseq{\I_{C_1/P_1}(dL_1)}{\I_{C_1/P_1}(dL_1+Q)}{I_{p_1,...,p_{e_2}/Q}(dL_1+Q)}.\]
Here the left term has vanishing $H^1$ by Maximal Rank, while the right term can be identified with
$(I_{p_1,...,p_{e_2}}/Q)((d-1)L_1)$ and
has vanishing $H^1$ by the generality of the points, in light of the bound $ e_2 \leq \binom{n-d+1}{d-1}$. 
Therefore the middle term has vanishing $H^1$
so the right map in \eqref{restriction} is surjective, hence so is the middle map.

This completes the proof of Theorem \ref{main-thm} in case $e$ is even. \par
The proof for 
$e=2e_2-1$ odd is
similar except we make $D_{e_2}$ a 'line' 
(i.e. a ruling of $P_1/\P^{n-1}$) instead of a conic by making $F_d$ go through 
the point $p_{e_2}$ as well. 
\end{proof}
\providecommand{\bysame}{\leavevmode\hbox to3em{\hrulefill}\thinspace}
\providecommand{\MR}{\relax\ifhmode\unskip\space\fi MR }
\providecommand{\MRhref}[2]{%
  \href{http://www.ams.org/mathscinet-getitem?mr=#1}{#2}
}
\providecommand{\href}[2]{#2}

\end{document}